\documentclass[12pt,oneside,english]{amsart}
\usepackage[latin1]{inputenc}
\usepackage{amssymb}

\makeatletter
\newtheorem{theorem}{Theorem}
\newtheorem{lemma}{Lemma}
\newtheorem{example}{Example}

\usepackage{babel}

\makeatother
\begin{document}
\baselineskip=17pt
\title[ Exponentially $S$-numbers]{Exponentially $S$-numbers}

\author{Vladimir Shevelev}
\address{Department of Mathematics \\Ben-Gurion University of the
 Negev\\Beer-Sheva 84105, Israel. e-mail:shevelev@bgu.ac.il}

\subjclass{11A51}

\begin{abstract}
Let $\mathbf{S}$  be the set of all finite or infinite increasing sequences
  of positive integers. For a sequence $S=\{s(n)\}, n\geq1,$ from $\mathbf{S},$ let
us call a positive number $N$ an exponentially $S$-number $(N\in E(S)),$
if all exponents in its prime power factorization are in $S.$ Let us accept that
$1\in E(S).$ We prove that, for every sequence $S\in \mathbf{S}$ with $s(1)=1,$
the exponentially $S$-numbers have a density
 $h=h(E(S))$ such that $$\sum_{i\leq x,\enskip i\in E(S)} 1 = h(E(S))x+O(\sqrt{x}\log x e^{c\frac{\sqrt{\log x}}{\log \log x}}),$$
where $c=4\sqrt{\frac{2.4}{\log 2}}=7.4430...$ and $h(E(S))=\prod_{p}(1+\sum_{i\geq2}\frac{u(i)-u(i-1)}{p^i}),$ where $u(n)$ is the characteristic
function of $S.$
\end{abstract}

\maketitle

\section{Introduction}
Let $\mathbf{S}$  be the set of all finite or infinite increasing sequences
of positive integers. For a sequence $S=\{s(n)\}, n\geq1,$ from $\mathbf{S},$ let
us call a positive number $N$ an exponentially $S$-number $(N\in E(S)),$
if all exponents in its prime power factorization are in $S.$ Let us accept that
$1\in E(S).$
For example, if $S=\{1\},$ then the exponentially 1-numbers form the sequence $B$ of
square-free numbers, and, as well-known,
\begin{equation}\label{1}
\sum_{i\leq x,\enskip i\in B} 1=\frac{6}{\pi^2}x+O(x^{\frac{1}{2}}).
\end{equation}
In case, when  $S=B,$ we obtain the exponentially square-free numbers (for the first time
this notion was introduced by M. V. Subbarao in 1972 \cite{6}, see A209061\cite{5}).
 Namely the exponentially square-free numbers were studied by many authors (for example, see
 \cite{2}, \cite{6} (Theorem 6.7), \cite{7}, \cite{8}, \cite{9}). In these papers,
 the authors analyzed the following asymptotic formula
 \begin{equation}\label{2}
\sum_{i\leq x,\enskip i\in E(B)} 1
=\prod_{p}(1+\sum_{a=4}^{\infty}\frac{\mu^2(a)-\mu^2(a-1)}{p^a})x+R(x),
\end{equation}
where the product is over all primes, $\mu$ is the M\"{o}bius function. The best
result of type $R(x)=o(x^{\frac{1}{4}})$ was obtained by Wu (1995) without using RH
(more exactly see \cite{9}). In 2007, assuming that RH is true, T\'{o}th \cite{8}
obtained $R(x)=O(x^{\frac{1}{5}+\varepsilon})$ and in 2010, Cao and Zhai \cite{2}
more exactly found that $R(x)=Cx^{\frac{1}{5}}+O(x^{\frac{38}{193}+\varepsilon}),$
where $C$ is a computable constant. Besides, T\'{o}th \cite{7} studied also the
exponentially $k$-free numbers, $k\geq2.$\newline
\indent In this paper, without using RH, we obtain a general formula with a
remainder term $O(\sqrt{x}\log x e^{c\frac{\sqrt{\log x}}{\log
\log x}})$ ($c$ is a constant) not depended on $S\in\mathbf{S}$ beginning with 1.
More exactly, we prove the following.
\begin{theorem}\label{t1}
For every sequence $S\in \mathbf{S}$ the exponentially $S$-numbers have a density
 $h=h(E(S))$ such that, $1)$ if $s(1)>1,$ then $h=0,$ while $2)$ if $s(1)=1,$ then
\begin{equation}\label{3}
\sum_{i\leq x,\enskip i\in E(S)} 1 = h(E(S))x+O(\sqrt{x}\log x e^{c\frac{\sqrt{\log x}}
{\log \log x}}),
\end{equation}
with $c=4\sqrt{\frac{2.4}{\log 2}}=7.443083...$ and
\begin{equation}\label{4}
h(E(S))=\prod_{p}(1+\sum_{i\geq2}\frac{u(i)-u(i-1)}{p^i}),
\end{equation}
where $u(n)$ is the characteristic function of sequence $S:\enskip u(n)=1,$ if\enskip
 $n\in S$ and $u(n)=0$ otherwise.
\end{theorem}
In particular, in case $S=B$ we obtain (\ref{2}) with a less good remainder term,
but which is suitable for all sequences in $\mathbf{S}$ beginning with 1.
\section{Lemma}
For proof Theorem \ref{t1} we need a lemma proved earlier (2007) by the author
\cite{4}, pp.200-202. For a fixed square-free number $r$, denote by $B_r$ the set of
square-free numbers $n$ for which $\gcd(n,r)=1$, and put
$$b_r(x)=|B_r\cap\{1,2,...,x\}|.$$
In particular, $B=B_1$ is the set of all square-free numbers.
\begin{lemma}\label{l1}
$$b_r(x)=\frac{6r}{\pi^2}\prod\limits_{p|r}(p+1)^{-1}x+ R_r(x),$$
\slshape where for every $x\geq1$ and every $r\in B$ \upshape

$$|R_r(x)|\leq \begin{cases}k\sqrt x,&\text{if \;$r\leq
N$}\\ke^{c\frac{\sqrt{\log r}}{\log\log
r}}\sqrt x,&\text{if\; $r\geq {N+1}$.}
\end{cases}$$

\slshape where\upshape\;\;\;$k=3.5\prod_{2\leq p\leq 23}(1+\frac{1}{\sqrt{p}})=
57.682607...$ (in case $r=1,\enskip k=3.5), \enskip c=4\sqrt{\frac{2.4}{\log 2}}=7.443083...,\enskip
N=6469693229.$

\end{lemma}
\section{Proof of Theorem 1}
1) Denote by $\Upsilon$ the sequence $\{2,3,4,...\}$ of all natural numbers without
 1. Let $S$ do not contain 1. Then, evidently, $E(S)\subseteq E(\Upsilon).$ Note that the
 sequence $E(\Upsilon)$ is called also powerful numbers (sequence A001694 in \cite{5}).
\newpage
 Bateman and Grosswald \cite{1} proved that

 \begin{equation}\label{5}
\sum_{i\leq x,\enskip i\in E(\Upsilon)} 1=\frac{\zeta(3/2)}{\zeta(3)} x^{1/2} +
\frac{\zeta(2/3)}{\zeta(2)} x^{1/3} + O(x^{1/6}).
\end{equation}
So, $h(E(\Upsilon))=0.$ Then what is more $h(E(S))=0.$ \newline\newline
Furthermore, \slshape denote by $r(n)$ the product of all distinct prime divisors of\upshape
\enskip $n;$ set
 $r(1)=1.$

\indent 2) Now let $1\in S.$ Note that the set $E(\Upsilon)\cap E(S)$ contains
1 and all numbers of $E(S)$ whose exponents in their prime power factorizations are more
than 1. Evidently, every number $y\in E(S)$ has
 a unique representation as the product of some
number $a\in E(\Upsilon)\cap E(S)$ and a number $m\in B_{r(a)}.$ In particular, if $y$ is square-free, then $a=1, \enskip m=y(\in B_1).$
For a fixed $a\in E(\Upsilon)\cap E(S),$ denote the set of $y=am\in E(S)$
by $E(S)^{(a)}.$ Then
$E(S)=\bigcup
\limits _{a\in E(S)\cap E(\Upsilon)} E(S)^{(a)}$, where the union is disjoint. Consequently,
by Lemma \ref{l1}, we have
\begin{equation}\label{6}
\sum_{i\leq x,\enskip i\in E(S)} 1=b_1 (x) + \sum_{4\leq a \leq x,\enskip
a\in E(S)\cap E(\Upsilon)} b_{r(a)}\left(\frac{x}{a}\right)
\end{equation}
$$ =\frac{6}{\pi^2} \left(1+\sum_{4\leq a \leq x,\enskip a\in E(S)\cap E(\Upsilon)}\;\prod_{p|r(a)}\left(1-\frac{1}{p+1}\right)\frac{1}
{a}\right)x+R(x),$$
where

\begin{equation*}
|R (x)|\leq 3.5\sqrt x + \sum_{4\leq a\leq x, \enskip a\in
E(S)\cap E(\Upsilon)}\left|R_{r(a)}\left (\frac {x}{a}\right)\right|
\leq 3.5\sqrt x
+
\end{equation*}

\begin{equation}\label{7}
+\sum_{\substack{4\leq a\leq x:r(a)\leq N\\a\in
E(S)\cap E(\Upsilon)}}\left|R_{r(a)}\left (\frac {x}{a}\right)\right|
+\sum_{\substack{a\leq x: r(a)\geq N+1\\a\in
E(S)\cap E(\Upsilon)}}\left|R_{r(a)}\left(\frac{x}{a} \right )\right |
\end{equation}

with  $N=6469693229$.

\indent Let $x>N$ go to infinity. Distinguish two cases:\enskip$(i) \enskip r(a)\leq N;
\enskip (ii)\enskip r(a)>N.$ \newline
(i) $r(a)\leq N.$ Denote by $E(\Upsilon)(n)$ the $n$-th powerful number (in increasing order).
According to (\ref{5}), $E(\Upsilon)(n)=(\frac{\zeta(3)}{\zeta(3/2)})^2n^2(1+o(1)).$
So, $\Sigma_{1\leq n\leq x}\frac{1}{\sqrt{E(\Upsilon)(n)}}=O(\log x).$ Hence,
by (\ref{7}) and Lemma \ref{l1},
$$ |R(x)|\leq3.5\sqrt{x}+k\sqrt{x}\sum_{a\leq x,\enskip a\in E(S)\cap E(\Upsilon)}\frac{1}{\sqrt{a}}=O(\sqrt{x}\log x).$$

(ii) $r(a)>N.$ Then, by (\ref{7}) and Lemma \ref{l1},
\newpage
$$R(x)\leq k\sqrt x\sum_{\substack{a\leq x :r(a)\geq N+1\\a\in E(S)\cap E(\Upsilon)}}
\frac{1}{\sqrt {a}} e^{c\frac{\sqrt{\log r(a)}}{\log\log r(a)}},$$
where the last sum does not exceed
$$
\sum\limits_{{N+1}\leq a\leq x:\enskip r(a)\geq {N+1}}\frac{1}{\sqrt{a}}
e^{c\frac{\sqrt{\log a}}{\log\log a}}\leq
e^{c\frac{\sqrt{\log x}}{\log\log x}}O(\log x).
$$

So, $R(x)=O(\sqrt{x}\log x e^{c\frac{\sqrt{\log x}}{\log \log x}})$
and, by (\ref{6}), we have

$$\sum_{i\leq x,\enskip i\in E(S)} 1=$$

$$\frac{6}{\pi^2} \left(1+\sum_{4\leq a \leq x,\enskip a\in E(S)\cap E(\Upsilon)}\;\prod_{p|r(a)}\left(1-\frac{1}{p+1}\right)\frac{1}
{a}\right)x+O(\sqrt{x}\log x e^{c\frac{\sqrt{\log x}}{\log \log x}}).$$

Moreover, if we replace here the sum $\sum_{a\leq x, a\in
E(S)\cap E(\Upsilon)}$ by the sum $\sum_{a\in E(S)\cap E(\Upsilon)},$ then the error
does not exceed $\frac {6x}{\pi^2}\sum\limits_{n>x}\frac{1}{E(\Upsilon)(n)}
=\frac {6x}{\pi^2} O(1/x)=O(1),$ then the result does not change. So, finally,
\begin{equation}\label{8}
\sum_{i\leq x,\enskip i\in E(S)} 1=
\end{equation}
$$\frac{6}{\pi^2} \left(\sum_{a\in E(S)\cap E(\Upsilon)}\;\prod_{p|r(a)}\left(1-\frac{1}{p+1}\right)\frac{1}
{a}\right)x+O(\sqrt{x}\log x e^{c\frac{\sqrt{\log x}}{\log \log x}}).$$
 Formula (\ref{8}) shows that, if $1\in S,$ then $E(S)$
has a density.
\section{Completion of the proof}
It remains to evaluate the sum (\ref{8}). For that we follow the scheme of \cite{4},
pp.203-204. For a fixed $l\in B$, denote by $C(l)$ the set of all
$E(S)\cap E(\Upsilon)$-numbers $a$ with $r(a)=l$. Recall that $r(1)=1.$ By
(\ref{8}), we have
\begin{equation}\label{9}
\sum_{i\leq x,\enskip i\in E(S)} 1=\frac{6}{\pi^2} x \sum\limits_{l\in B} \prod\limits_
{p|l}\left(1-\frac{1}{p+1}\right)\sum\limits_{a\in
C(l)}\frac{1}{a}+R(x).
\end{equation}

Consider the function $A:N\rightarrow R$ given by:
\begin{equation*}
A(l)=\begin{cases}
\sum\limits_{a\in
C(l)}\frac{1}{a},&l\in B,\\
0,&l\not\in B.
\end{cases}
\end{equation*}
\begin{example}\label{e1}
$$A(1)=\sum\limits_{a\in C(1)}\frac{1}{a}=\sum\limits_{r(a)=1}
\frac{1}{a}=1.$$
\end{example}
\newpage
\begin{example}\label{e2}
Let $p$ be prime. Since $r(p)=p,$ then
$$A(p)=\sum\limits_{a\in C(p)}\frac{1}{a}=\sum\limits_{i\geq2}\frac{1}{p^{s(i)}}.$$
The sum not contains $\frac{1}{p^{s(1)}}=\frac{1}{p}$ since, by the condition,
$a\in E(S)\cap E(\Upsilon),$ but the sequence $E(\Upsilon)$ not contains any prime.
\end{example}
\begin{example}\label{e3}
Let $p<q$ be primes. Since $r(pq)=pq,$ then
$$A(pq)=\sum\limits_{i\geq2, j\geq2}\frac{1}{p^{s(i)}}\frac{1}{q^{s(j)}}. $$
\end{example}
It is evident that, if $l_1,l_2\in B$ and $\gcd(l_1,l_2)=1$, then
\begin{equation*}
A(l_1l_2)=\sum\limits_{a\in
C(l_1l_2)}\frac{1}{a}=\sum\limits_{a\in
C(l_1)}\frac{1}{a}\sum\limits_{a\in
C(l_2)}\frac{1}{a}=A(l_1)A(l_2).
\end{equation*}

It follows that $A(l)$ is a multiplicative function. Hence the
function $f$ which is defined by

\begin{equation*}
f(l)=\prod\limits_{p|l} \left(1-\frac{1}{p+1}\right) A(l)
\end{equation*}

is also multiplicative. Evidently, by the definition of $A(n),$
$$\sum_{n=1}^\infty f(n)\leq\sum_{n=1}^\infty A(n)\leq\sum_{a\in E(\Upsilon)} \frac{1}{a}<\infty. $$
Consequently (\cite{3},p.103):
\begin{equation}\label{10}
\sum_{n=1}^\infty f(n)=\prod\limits_{p}(1+f(p)+f(p^2)+\ldots).
\end{equation}

Since $f(p^k)=0$ for $k\geq 2$, then by (\ref{9}):

\begin{equation*}
\sum_{i\leq x,\enskip i\in E(S)} 1=\frac{6}{\pi^2} x\sum_{l=1}^\infty f(l)+R(x)
=\frac{6}{\pi^2}
x\prod_{p}(1+f(p))+R(x)=
\end{equation*}

\begin{equation*}
=\frac{6}{\pi^2} x\prod_{p
}\left(1+\left(1-\frac{1}{p+1}\right)\left(\frac{1}{p^{s(2)}}+
\frac{1}{p^{s(3)}}+\frac{1}{p^{s(4)}}+\ldots\right)\right)+R(x).
\end{equation*}

Now we have
\begin{equation}\label{11}
h(E(S))=\frac{6}{\pi^2}\prod_{p}(1+(1-\frac{1}{p+1})\sum_{i\geq2} \frac{1}
{p^{s(i)}})=
\end{equation}
$$\frac{6}{\pi^2}\prod_{p}(1+\sum_{i\geq2} \frac{p}
{(p+1)p^{s(i)}})=$$
\newpage
$$\prod_{p}((1-\frac{1}{p^2})-(1-\frac{1}{p})\sum_{i\geq2}\frac{u(i)}{p^i})$$
and, taking into account that $u(1)=1,$ we find
$$h(E(S))=\prod_{p}(1-\frac{1}{p^2}-(1-\frac{1}{p})\frac{1}{p}+(1-\frac{1}{p})\sum_{i\geq1}
\frac{u(i)}{p^i})=$$
$$\prod_{p} ( (1-\frac{1}{p})+
\sum_{i\geq1}\frac{u(i)}{p^{j}} -\frac{1}{p}\sum_{i\geq1}\frac{u(i)}{p^{j}})=$$
$$\prod_{p}((1-\frac{1}{p})+\frac{1}{p}+
\sum_{i\geq2}\frac{u(i)}{p^{j}} -\frac{1}{p}\sum_{i\geq2}\frac{u(i-1)}{p^{i-1}})= $$
$$\prod_{p}(1+\sum_{i\geq2}\frac{u(i)-u(i-1)}{p^i})$$
 which gives the required evaluation of the sum in (\ref{8}) and completes the
 proof of the theorem.

\section{A question of D. Berend}

Let $p_n$ be the $n$-th prime. Let $\mathrm{A}=\{S_1,S_2,...\}$ be an infinite
sequence of sequences  $S_i\in\mathbf{S}$ beginning with 1. We say that a positive
number $N$ is an exponentially $\mathrm{A}$-number $(N\in E(\mathrm{A})),$
if in case that $p_n,\enskip n\geq1,$
divides $N,$ then its exponent in the prime power factorization of $N$ belongs
to $S_n.$ We accept that $1\in E(\mathrm{A}).$ How will change Theorem \ref{t1}
for the exponentially $\mathrm{A}$-numbers?
\newline
\indent An analysis of the proof of Theorem \ref{t1} shows that also in this more general case,
for every sequence $\mathrm{A}$ there exists a density $h(\mathrm{A})$ of the
exponentially $\mathrm{A}$-numbers such that
 \begin{equation}\label{12}
\sum_{i\leq x,\enskip i\in E(\mathrm{A})} 1 = h(E(\mathrm{A}))x+R(x),
\end{equation}
where $R(x)$ is the same as in Theorem \ref{t1} and
\begin{equation}\label{13}
h(E(\mathrm{A}))=\prod_{n\geq1}(1+\sum_{i\geq2}\frac{u_n(i)-u_n(i-1)}{{p_n}^i}),
\end{equation}
where $u_n(k)$ is the characteristic function of sequence $S_n:\enskip u_n(k)=1,$ if
\enskip $k\in S_n$ and $u_n(k)=0$ otherwise.
\begin{example}\label{e4}
Let
$$ \mathrm{A}=\{S_1=\{1\}, S_2=\{1,2\}, ... , S_n=\{1,...,n\}, ... \}.$$
Then, by $(\ref{13}),$
$$h(E(A))=\prod_{n\geq1}(1-\frac{1}{{p_n}^{n+1}})=0.7210233...\enskip .$$
\end{example}
\newpage
\section{A question}

\indent Let $1\in S.$ Then the density $h(E(S))$ is in the interval $[6/\pi^2, 1].$
Whether the set $\{h(E(S))\}$ is a dense set in this interval?\newline
\indent D. Berend (private communication) gave a negative answer. Indeed, consider the set
$\mathbf{S}_1$ of sequences $\{S\}$ containing 2. Then, evidently,
$h(E(S))\geq h(E(\{1,2\}))$ such that, by Theorem \ref{t1},

\begin{equation}\label{14}
h(E(S))|_{S\in \mathrm{S_1}}\geq\prod_{p}(1-\frac{1}{p^3}).
\end{equation}
Now consider  the set $\mathbf{S}_2$ of sequences $\{S\}$ not containing 2. Then
$h(E(S))\leq h(E(\{1,3,4,5,6,...\}))$ such that, by Theorem \ref{t1},
\begin{equation}\label{15}
h(E(S))|_{S\in \mathrm{S_2}}\leq\prod_{p}(1-\frac{1}{p^2}+\frac{1}{p^3})=
\prod_{p}(1-\frac{p-1}{p^3}).
\end{equation}

Thus, by (\ref{14})-(\ref{15}), we have a gap in the set $\{h(E(S))\}$
in interval
$$(\prod_{p}(1-\frac{p-1}{p^3}),\enskip \prod_{p}(1-\frac{1}{p^3}) ). $$
\newline\newline
\indent Of course, this Berend's idea has far-reaching effects.

\section{Acknowledgement}
The author is very grateful to Jean-Paul Allouche for sending the important papers
\cite{6}, \cite{7} and \cite{9}, and to Daniel Berend for the question in Section 5
and the answer on author's question in Section 6. \newline

\newpage
\end{document}